\documentclass[14pt,twoside]{article}
\usepackage{graphicx}
\usepackage{amsmath,amsbsy,amsfonts}
\usepackage{color}
\def\dis{\displaystyle}
\def\nd{\noindent}

\newtheorem{lemma}{Lemma}[section]
\newtheorem{prop}{Proposition}[section]
\newtheorem{theorem}{Theorem}[section]

\newtheorem{coro}{Corollary}[section]

\newtheorem{definition}[theorem]{Definition}

\newcommand{\ds}{\displaystyle}
\def\dis{\displaystyle}

\let\Section=\section
\def\section{\setcounter{equation}{0}\Section}

\begin{document}
\title{ On a $\Phi$-Kirchhoff multivalued problem with critical growth in an Orlicz-Sobolev space}

\author{
Giovany M. Figueiredo\thanks{Partially supported by CNPq/PQ
301242/2011-9 and 200237/2012-8 }  \\
\noindent Universidade Federal do Par\'a, Faculdade de Matem\'atica, \\
\noindent CEP: 66075-110, Bel\'em - Pa, Brazil \\
\noindent e-mail: giovany@ufpa.br  \vspace{0.5cm}\\
\noindent
Jefferson A. Santos \thanks{Partially supported by CNPq-Brazil grant Casadinho/Procad 552.464/2011-2 } \\
\noindent Universidade Federal de Campina Grande,\\
\noindent Unidade Acad\^emica de Matem\'atica e Estat\'istica,\\
\noindent CEP:58109-970, Campina Grande - PB, Brazil\\
\noindent e-mail: jefferson@dme.ufcg.edu.br\vspace{0.5cm}\\
}\vspace{0.5cm}
\date{}

\pretolerance10000

\maketitle

\begin{abstract}
{This paper is concerned with the multiplicity of nontrivial
solutions in an Orlicz-Sobolev space for a nonlocal problem with critical growth,
involving N-functions and theory of locally Lispchitz continuous
functionals. More precisely, in this paper, we study a result of
multiplicity to the following multivalued elliptic problem:
$$
\left \{ \begin{array}{l}
-M\left(\displaystyle\int_\Omega \Phi(|\nabla u|)dx\right)\Delta_\Phi u\in  \partial F(.,u)+\alpha h(u)  \ \mbox{in}\ \Omega,\\
u\in W_0^1L_\Phi(\Omega),                                    
\end{array}
\right.
$$
where  $\Omega\subset\mathbb{R}^{N}$ is a bounded smooth domain,
$N\geq 3$, $M$ is a continuous function, $\Phi$ is an  N-functions,
$h$  is an odd increasing homeomorphism from $\mathbb{R}$ to
$\mathbb{R}$, $\alpha$ is positive parameter, $\Delta_\Phi$ is the
corresponding $\Phi-$Laplacian and $\partial F(.,t)$ stands for
Clarke generalized of a function $F$ linked  with critical growth.
We use genus theory to obtain the main result.}
\end{abstract}

\maketitle

\section{Introduction}

Let $\Omega\subset\mathbb{R}^{N}$, $N\geq 3$, be a bounded
 domain with smooth boundary $\partial \Omega$ and consider a continuous
 function $\phi:(0,+\infty)\rightarrow [0,+\infty)$.
The purpose of this article is to investigate the multiplicity of
nontrivial solutions to the multivalued elliptic problem
$$
\left \{ \begin{array}{l}
-M\left(\displaystyle\int_\Omega \Phi(|\nabla u|)dx\right)\Delta_\Phi u\in  \partial F(.,u) +\alpha h(u) \ \mbox{in}\ \Omega,\\
u\in W_0^1L_\Phi(\Omega),                                    
\end{array}
\right.\leqno{(P_\alpha)}
$$
 \noindent where $\alpha>0$ is a parameter, $\Phi$ is defined for $t\in \mathbb{R}$ by setting
$$
\Phi(t)=\int_0^{|t|}\phi(s)sds,
$$
 $\Delta_\Phi$ is the operator,
$$
\Delta_\Phi u:=div(\phi(|\nabla u|)\nabla u),
$$
$F:\Omega\times\mathbb{R}\rightarrow\mathbb{R}$ is mensurable and
$$
\partial F(x,t)=\left\{s\in\mathbb{R}; F^0(x,t;r)\geq sr, \ r\in \mathbb{R}\right\}.
$$
Here $F^0(x,t;r)$ denotes the generalized directional derivative of
$t\mapsto F(x,t)$ in direction of $r$, that is,
$$
F^0(x,t;r)=\displaystyle \limsup_{h\rightarrow
t,s\downarrow0}\frac{F(x,h+sr)-F(x,h)}{s}.
$$
In this paper we consider
$$
F(x,t)=\int_0^tf(x,s)ds,
$$
where $f:\Omega\times \mathbb{R}\rightarrow \mathbb{R}$ is
mensurable and $f(x,.)$ is locally bounded in $\mathbb{R}$ and
$$
\underline{f}(x,t)=\displaystyle \lim_{\epsilon\downarrow 0}\mbox{ess
inf}\left\{ f(x,s);|s-t|<\epsilon\right\} \mbox{ and }
\overline{f}(x,t)=\lim_{\epsilon\downarrow 0}\mbox{ess sup}\left\{
f(x,s);|s-t|<\epsilon\right\} .
$$
It is well known that
\begin{equation}\label{partial1}
\partial F(x,t)=[\underline{f}(x,t),\overline{f}(x,t)], \mbox{ (see
\cite{chang})},
\end{equation}
and that, if $f(x,.)$ is continuous then $\partial F(x,t)=\{f(x,t)\}$.

Problem $(P_\alpha)$  with $N=3$ and $\phi(t)=4$, that is,
$$
\left \{ \begin{array}{l}
-M\left(\displaystyle\int_\Omega |\nabla u|^{2} dx\right)\Delta u \in \partial F(.,u)+\alpha h(u) \ \mbox{in}\ \Omega,\\
u\in H^{1}_{0}(\Omega)                                    
\end{array}
\right.\leqno{(*)}
$$
is called nonlocal because of the presence of the term
$M\left(\dis\int_{\Omega}|\nabla u|^{2} dx \right)$ which implies
that the equation $(*)$ is no longer a pointwise identity. This
phenomenon causes some mathematical difficulties which makes the
study of such a class of problem particularly interesting. Besides,
this class of problems has physical motivation. Indeed, the operator
$M(\int_\Omega|\nabla u|^{2} dx )\Delta u$ appears in the Kirchhoff
equation, which arises in nonlinear vibrations, namely
$$
\left\{
\begin{array}{l}
 u_{tt}-M\left(\dis\int_{\Omega}|\nabla u|^{2} dx \right)\Delta u =
g(x,u) \ \mbox{in} \ \Omega \times (0,T)\\
u=0 \ \mbox{on} \ \partial\Omega \times (0,T)\\
u(x,0)=u_{0}(x)\ \ , \ \ u_{t}(x,0)=u_{1}(x).
\end{array}
\right.\leqno{(1.1)}
$$

Such a hyperbolic equation is a general version of the Kirchhoff
equation
$$
\rho\displaystyle\frac{\partial^{2}u}{\partial
t^{2}}-\biggl(\displaystyle\frac{P_{0}}{h}+\displaystyle\frac{E}{2L}\int^{L}_{0}\biggl|\displaystyle\frac{\partial
u}{\partial x}\biggl|^{2} dx
\biggl)\displaystyle\frac{\partial^{2}u}{\partial x^{2}}=0
\leqno{(1.2)}
$$
presented by Kirchhoff \cite{kirchhoff}. This equation extends the
classical d'Alembert's wave equation by considering the effects of
the changes in the length of the strings during the vibrations. The
parameters in equation $(1.2)$ have the following meanings: $L$ is
the length of the string, $h$ is the area of cross-section, $E$ is
the Young modulus of the material, $\rho$ is the mass density and
$P_{0}$ is the initial tension.

When an elastic string with fixed ends is subjected to transverse
vibrations, its length varies with the time: this introduces changes
of the tension in the string. This induced Kirchhoff to propose a
nonlinear correction of the classical D'Alembert equation. Later on,
Woinowsky-Krieger (Nash - Modeer) incorporated this correction in
the classical Euler-Bernoulli equation for the beam (plate) with
hinged ends. See, for example, \cite{Arosio}, \cite{Arosio1} and the
references therein.

This class of problems began to call attention of several
researchers mainly after the work of Lions $\cite{lions}$, where a
functional analysis approach was proposed to attack it.

We have to point out that nonlocal problems also appear in other
fields as, for example, biological systems where $u$ describes a
process which depends on the average of itself (for example,
population density). See, for example, \cite{alvescorreama} and its
references.

The reader may consult $\cite{alvescorreama}$, $\cite{alvescorrea}$,
 $\cite{GJ}$, $\cite{ma}$ and the references therein, for more
information on nonlocal problems.

On the other hand, in this study, the nonlinearity $f$ can be
discontinuous. The interest in the study of nonlinear partial
differential equations with discontinuous nonlinearities has
increased because many free boundary problems arising in
mathematical physics may be stated in this form. Among these
problems, we have the obstacle problem, the seepage surface problem,
and the Elenbaas equation, see for example \cite{chang},
\cite{chang1} and \cite{chang2}.

For enunciate the main result, we need to give some hypotheses on
the functions  $M, \phi$, $h$ and $f$.

The hypotheses on the function $\phi:(0,+\infty)\rightarrow
(0,+\infty)$ of $C^{1}$ class are the following:

\begin{description}

\item[($\phi_1$)] For all $t>0$,
$$
\phi(t)>0 \ \ \mbox{and} \ \ (\phi(t)t)'>0.
$$

\item[($\phi_2$)] There exist $l \in (\frac{N}{2},N)$,  $l< m < \frac{l^{*}}{2}=
\displaystyle\frac{lN}{2(N-l)}$ such that
$$
l\leq \frac{\phi(t)t^{2}}{\Phi(t)}\leq m,
$$
for $t> 0$, where $\Phi(t)=\displaystyle\int^{|t|}_{0}\phi(s)s ds$.
\end{description}
Through this paper, we assume that $h$  is an odd increasing
homeomorphism from $\mathbb{R}$ to $\mathbb{R}$, verifying
\begin{description}
\item[($H_0$)] There exist $1<h_0<h_1<l$ such that
$$
h_0\leq \frac{h(t)t}{H(t)}\leq h_1, \ t>0,
$$
where $H(t)=\int_0^th(s)ds$.
\end{description}
Considering $\widetilde{H}(t):=\int_0^th^{-1}(s)ds$, then we can
obtain complementary N-function which define corresponding Orlicz
space $L_H(\Omega)$.

The hypothesis on the continuous function
$M:\mathbb{R}^{+}\rightarrow \mathbb{R}^{+}$ is the following:

\begin{description}
\item[($M_1$)]
There exists $\sigma>0$ such that
$$
M(t)\geq \sigma \mbox{ for all }t\geq 0.
$$
\item[($M_2$)] There exists $\theta\in(2m,l^*)$ such that
$$
\widehat{M}(t)-\frac{m}{\theta} M(t)t\geq 0, \ \mbox{ for all }t\geq 0,
$$
where $\widehat{M}(t)=\int_0^tM(s)ds.$

\end{description}

The hypotheses on the function $f:\Omega\times\mathbb{R}\rightarrow \mathbb{R}$
are the following:
\begin{description}

\item[($f_1$)] $f(x,t)$ is odd with respect $t$, there exist $b_{0},b_{1}>0$ and $a_0\geq 0$ such that
$$
b_0\phi_*(t)t\leq f(x,t)\leq b_1\phi_*(t)t,\ t\geq a_0 \text{ and }x\in \Omega,
$$
with $1<\frac{b_1}{b_0}<\frac{l^*}{\theta}$, where $\phi_*(t)t$ is such that Sobolev conjugate function $\Phi_*$ of $\Phi$ is its primitive, that is, $\Phi_*(t)=\int_0^{|t|}\phi_*(s)sds$.

\item[($f_2$)] Moreover, if $a_0>0$, we consider $f(x,t)=0$ for all  $|t|< a_0$ and $x\in \Omega$.
\end{description}

The main result of this paper is:

\begin{theorem}\label{teorema1}
Assume that conditions $(\phi_{1})$, $(\phi_{2})$, $(H_0)$, $(M_{1})$, $(M_{2})$,
$(f_{1}),(f_{2})$ hold. If $\underline{f}$ and $\overline{f}$ are N-functions, then for $\alpha=\alpha(\theta,h_0,b_0,b_1,l,m, \sigma, S_N)>0$, $a_0=a_0(H(1),\sigma,l,|\Omega|, \alpha, h_1)>0$ sufficiently small (or $a_0=0$), the problem $(P_\alpha)$  has infinitely many solutions.
\end{theorem}
{\bf An Application: } Consider the equation
$$
\left \{ \begin{array}{l}
-\left(a\displaystyle\int_\Omega \Phi(|\nabla u|)dx+b\right)\Delta_\Phi u=\widehat{H}(|u|-a_0) \phi_*(|u|)u+\alpha h(u)  \ \mbox{ in}\ \Omega,\\
u\in W_0^1L_\Phi(\Omega),
\end{array}
\right. \leqno{(A_1)}
$$
where $a,b,a_0>0$ and $\widehat{H}$ is the Heaviside function. We claim
that  ($A_1$)
 admits infinitely many solutions in $W_0^{1,\Phi}(\Omega)$, for $\alpha,a_0>0$ sufficiently small. Indeed, let $f(x,t)=\widehat{H}(|t|-a_0) \phi_*(|t|)t$ and set
$$
F(x,u)=\int_0^u f(x,t)dt=\chi_{[|u|\geq a_0]}\left(\Phi_*(u)-\Phi_*(a_0)\right).
$$
Next consider the associated multivalued equation
$$
\left \{ \begin{array}{l}
-\left(a\displaystyle\int_\Omega \Phi(|\nabla u|)dx+b\right)\Delta_\Phi u\in\partial F(.,u)+\alpha h(u)  \ \mbox{ in}\ \Omega.
\end{array}
\right. \leqno{(A_2)}
$$
We will show that the assumptions of Theorem \ref{teorema1} hold
true. Indeed, since $\Omega$ is bounded, it is straightforward to
check that conditions $(M_1),(M_2),(f_1)$ and $(f_2)$ hold. Thus by Theorem
\ref{teorema1}, ($A_2$) admits infinitely many solutions in
$W_0^{1,\Phi}(\Omega)$ with $\alpha,a_0>0$ sufficiently small.
Finally, since that $|[|u|=a_0]|=0$ the problem ($A_1$)  has
infinitely many solutions.

In the last twenty years the study on nonlocal problems of the type
$$
\left \{ \begin{array}{l}
-M\left(\displaystyle\int_\Omega \mid\nabla u\mid^{2} dx\right)\Delta u= g(x,u)  \ \mbox{in}\ \Omega,\\
u\in H^{1}_{0}(\Omega)                                    
\end{array}
\right.\leqno{(K)}
$$
grew exponentially. That was, probably, by the difficulties existing
in this class of problems and that do not appear in the study of
local problems, as well as due to their significance in
applications. Without hope of being thorough, we mention some
articles with multiplicity results and that are related with our
main result. We will restrict our comments to the works that have
emerged in the last four years

The problem $(K)$ was studied in \cite{GJ}. The version with
p-Laplacian operator was studied in \cite{correa}. In both cases,
the authors showed a multiplicity result using genus theory. In
\cite{Xiao} and \cite{Ji} the authors showed a multiplicity result
for the problem $(K)$ using the Fountain theorem and the Symmetric
Mountain Pass theorem. In all these articles the nonlinearity is
continuous. The case discontinuous was studied in \cite{Nascimento}.
With a nonlinearity of the Heaviside type the authors showed a
existence of two solutions via Mountain Pass Theorem and Ekeland's
Variational Principle.

In \cite{Chung1} the author showed the existence of two solutions
for a problem involving $\Phi$-Kirchhoff operator and nonlinearity
of variable exponent type and subcritical growth. He used the
mountain pass theorem combined with the minimum principle.

In this work we extend the studies found in the papers above in the
following sense:

\noindent a) We cannot use the classical Clark's Theorem for $C^1$
functional (see \cite[Theorem 3.6]{Davi}), because in our case, the
energy functional is only locally Lipschitz continuous. Thus, in all
section \ref{finalissimo} we adapt for nondifferentiable functionals
an argument found in \cite{GP}.\\

\noindent b) Unlike \cite{Chung1} and \cite{Nascimento}, we show a
result of multiplicity using genus theory considering a nonlinearity
that can have a number enumerable of discontinuities.\\

\noindent c) Problem $(P_\alpha)$ possesses more complicated
nonlinearities, for example:

\noindent (i) $\Phi(t)=t^{p_{0}}+ t^{p_{1}}$, $1< p_{0} < p_{1}< N$
and $ p_1\in (p_0,p^{*}_{0})$.

\noindent (ii) $\Phi(t)=(1+t^{2})^{\gamma}-1$, $\gamma\in
(1,\frac{N}{N-2})$.

\noindent (iii) $\Phi(t)=t^{p}\log(1+t)$ with $1<p_0<p<N-1$, where
$p_0=\frac{-1+\sqrt{1+4N}}{2}$.

\noindent (iv) $\Phi(t)=\int_0^ts^{1-\theta}\left(\sinh^{-1}s\right)^\beta ds, \ 0\leq\theta\leq 1, \ \beta>0$.\\

Before concluding this introduction, it is very important to say
that in the literature, we find many papers where the authors study
the existence and multiplicity of solution for problems involving
the $\Phi-$Laplacian operator, see, for example, \cite{Bonanno},
\cite{Carvalho}, \cite{Chung}, \cite{Mihailescu}, \cite{Santos},
and references therein.

The paper is organized as follows. In the next section we present a
brief review on Orlicz-Sobolev spaces. In section  \ref{Section
results for the discontinuous} we recall some definitions  and basic
results on the critical point theory of locally Lipschitz continuous
functionals. In Section \ref{Section Genus} we present just some preliminary results
involving genus theory that will be used in this work. In the
Section \ref{finalissimo} we prove Theorem \ref{teorema1}.


\section{A brief review on Orlicz-Sobolev spaces}\label{Section Orlicz}

Let $\varphi$ be a real-valued function defined in $[0,\infty)$ and
having the following properties:

\noindent $a)$ \ \  $\varphi(0)=0$, $\varphi(t)>0$ if $t>0$ and
$\displaystyle\lim_{t\rightarrow \infty}\varphi(t)=\infty$.

\noindent $b)$ \ \  $\varphi$ is nondecreasing, that is, $s>t$ implies
$\varphi(s) \geq \varphi(t)$.

\noindent $c)$ \ \ $\varphi$ is right continuous, that is,
$\displaystyle\lim_{s\rightarrow t^{+}}\varphi(s)=\varphi(t)$.

Then, the real-valued function $\Phi$ defined on $\mathbb{R}$ by
$$
\Phi(t)= \displaystyle\int^{|t|}_{0}\varphi(s) \ ds
$$
is called an N-function. For an N-function $\Phi$ and  an open set
$\Omega \subseteq \mathbb{R}^{N}$, the Orlicz space
$L_{\Phi}(\Omega)$ is defined (see \cite{adams}). When $\Phi$
satisfies $\Delta_{2}$-condition, that is, when there are $t_{0}\geq
0$ and $K>0$ such that $\Phi(2t)\leq K\Phi(t)$, for all $t\geq
t_{0}$,  the space $L_{\Phi}(\Omega)$ is the vectorial space of the
measurable functions $u: \Omega \to \mathbb{R}$ such that
$$
\displaystyle\int_{\Omega}\Phi(|u|) \ dx < \infty.
$$
The space $L_{\Phi}(\Omega)$ endowed with Luxemburg norm, that is,
the norm given by
$$
|u|_{\Phi}= \inf \biggl\{\tau >0:
\int_{\Omega}\Phi\Big(\frac{|u|}{\tau}\Big)\ dx\leq 1\biggl\},
$$
is a Banach space. The complement function of $\Phi$, denoted by
$\widetilde{\Phi}$, is given by the Legendre transformation, that is
$$
\widetilde{\Phi}(s)=\displaystyle\max_{t \geq 0}\{st -\Phi(t)\} \ \
\mbox{for} \ \ s \geq 0.
$$
These $\Phi$ and $\widetilde{\Phi}$ are complementary each other.
Involving the functions $\Phi$ and $\widetilde{\Phi}$, we have the
Young's inequality given by
$$
st \leq \Phi(t) + \widetilde{\Phi}(s).
$$
Using the above inequality, it is possible to prove the following
H\"{o}lder type inequality
$$
\biggl|\displaystyle\int_{\Omega}u v \ dx \biggl|\leq
2|u|_{\Phi}|v|_{\widetilde{\Phi}}\,\,\, \forall \,\, u \in
L_{\Phi}(\Omega) \,\,\, \mbox{and} \,\,\,  v \in
L_{\widetilde{\Phi}}(\Omega).
$$

Hereafter, we denote by $W^{1}_{0}L_{\Phi}(\Omega)$ the
Orlicz-Sobolev space obtained by the completion of
$C^{\infty}_{0}(\Omega)$ with norm
$$
\|u\|_\Phi=|u|_\Phi+|\nabla u|_\Phi.
$$

When $\Omega$ is bounded, there is $c>0$ such that
$$
|u|_\Phi \leq c |\nabla u|_\Phi.
$$

In this case, we can consider

$$
\|u\|_{\Phi}=|\nabla u|_{\Phi}.
$$

Another important function related to function $\Phi$, is the
Sobolev conjugate function $\Phi_{*}$ of $\Phi$ defined by
$$
\Phi^{-1}_{*}(t)=\displaystyle\int^{t}_{0}\displaystyle\frac{\Phi^{-1}(s)}{s^{(N+1)/N}}ds, \ t>0.
$$

Another important inequality  was proved by    Donaldson and Trudinger \cite{Donaldson2}, which establishes that
for all open $\Omega \subset \mathbb{R}^{N}$ and there is a constant  $S_N=S(N) > 0$ such that
\begin{equation}\label{trudinger-emb}
\mid u\mid_{\Phi_*}\leq S_N\mid\nabla
u\mid_{\Phi}, \ u\in W_0^{1,\Phi}(\Omega).
\end{equation}
This inequality shows the below embedding is continuous
$$
W_0^{1,\Phi}(\Omega) \stackrel{\hookrightarrow}{\mbox{\tiny cont}} L_{\Phi_*}(\Omega).
$$
If bounded domain $\Omega$ and the limits below hold
\begin{equation} \label{M1}
\limsup_{t \to 0}\frac{B(t)}{\Phi(t)}< +\infty \,\,\, \mbox{and} \,\,\, \limsup_{|t| \to +\infty}\frac{B(t)}{\Phi_{*}(t)}=0,
\end{equation}
the embedding
\begin{equation} \label{M2}
W_0^{1,\Phi}(\Omega) \hookrightarrow L_{B}(\Omega)
\end{equation}
is compact.

The hypotheses $(\phi_{1})-(\phi_{2})$ implies that $\Phi$,
$\widetilde{\Phi}$, $\Phi_{*}$ and $\widetilde{\Phi}_{*}$ satisfy
$\Delta_{2}$-condition. This condition allows us conclude that:

\noindent 1) $u_{n}\rightarrow 0$ in $L_{\Phi}(\Omega)$ if, and only
if, $\displaystyle\int_{\Omega}\Phi(u_{n})\ dx\rightarrow 0$.

\noindent 2)  $L_{\Phi}(\Omega)$ is separable and
$\overline{C^{\infty}_{0}(\Omega)}^{|.|_{\Phi}}=L_{\Phi}(\Omega)$.

\noindent 3)  $L_{\Phi}(\Omega)$ is reflexive and its dual is
$L_{\widetilde{\Phi}}(\Omega)$(see \cite{adams}).

\vspace{.5cm}

Under assumptions $(\phi_{1})-(\phi_{2})$, some elementary
inequalities listed in the following lemmas are valid. For the
proofs, see \cite{fukagai}.

\begin{lemma}
Assume $(\phi_1)-(\phi_2)$. Then,
$$
\Phi(t) = \int_0^{|t|} s \phi(s) ds, \,\,\,
$$
is a $N$-function with $\Phi, \widetilde{\Phi} \in \Delta_2$.  Hence, $L_\Phi(\Omega), W^{1,\Phi}(\Omega)$ and $W_0^{1,\Phi}(\Omega)$ are reflexive and separable spaces.
\end{lemma}

\begin{lemma}\label{DESIGUALD}The functions $\Phi$, $\Phi_*$,  $\widetilde{\Phi}$ and $\widetilde{\Phi}_*$ satisfy the inequality
\begin{equation} \label{D1}
\widetilde{\Phi}(\phi(t)t) \leq \Phi(2t)  \mbox{ and } \widetilde{\Phi}_*(\phi_*(t)t) \leq \Phi_*(2t)\,\,\, \forall t \geq 0.
\end{equation}
\end{lemma}

\begin{lemma} \label{desigualdadeimportantes} Assume that $(\phi_1)-(\phi_2)$ hold and let $\xi_{0}(t)=\min\{t^{l},t^{m}\}$,\linebreak $ \xi_{1}(t)=\max\{t^{l},t^{m}\},$ for all $t\geq 0$. Then,
$$
\xi_{0}(\rho)\Phi(t) \leq \Phi(\rho t) \leq  \xi_{1}(\rho)\Phi(t) \;\;\; \mbox{for} \;\; \rho, t \geq 0
$$
and
$$
\xi_{0}(|u|_{\Phi}) \leq \int_{\Omega}\Phi(u) \leq \xi_{1}(|u|_{\Phi})  \;\;\; \mbox{for} \;\; u \in L_{\Phi}(\Omega).
$$
\end{lemma}

\begin{lemma} \label{F3} The function $\Phi_*$ satisfies the following inequality
$$
l^{*} \leq \frac{\Phi'_*(t)t}{\Phi_{*}(t)} \leq m^{*} \,\,\, \mbox{for} \,\,\, t > 0.
$$
\end{lemma}
As an immediate consequence of the Lemma \ref{F3}, we have the following result

\begin{lemma} \label{F2} Assume that $(\phi_1)-(\phi_2)$ hold and let  $\xi_{2}(t)=\min\{t^{l^{*}},t^{m^{*}}\},$ $\xi_{3}(t)=\max\{t^{l^{*}},t^{m^{*}}\}$ for all $t\geq 0$. Then,
$$
\xi_{2}(\rho)\Phi_*(t) \leq \Phi_*(\rho t) \leq \xi_{3}(\rho)\Phi_*(t) \;\;\; \mbox{for} \;\; \rho, t \geq 0
$$
and
$$
\xi_{2}(|u|_{\Phi_*}) \leq \int_{\Omega}\Phi_*(u)dx \leq \xi_{3}(|u|_{\Phi_*})  \;\;\; \mbox{for} \;\; u \in L_{A_*}(\Omega).
$$
\end{lemma}
\begin{lemma} \label{lem Phiest}
Let $\widetilde{\Phi}$ be the complement of $\Phi$ and put
$$
\xi_4(s)=\min\{s^{\frac{l}{l-1}}, s^{\frac{m}{m-1}}\}\ \mbox{and}\
\xi_5(s)=\max\{s^{\frac{l}{l-1}}, s^{\frac{m}{m-1}}\}, \ s\geq0.
$$
Then the following inequalities hold
$$
\xi_4(r)\widetilde{\Phi}(s)\leq\widetilde{\Phi}(r s)\leq
\xi_5(r)\widetilde{\Phi}(s),\ r,s\geq0
$$
and
$$
\xi_4(| u|_{\widetilde{\Phi}})\leq
\int_\Omega\widetilde{\Phi}(u)dx\leq \xi_5(|
u|_{\widetilde{\Phi}}),\ u\in
L_{\widetilde{\Phi}}(\Omega).
$$
\end{lemma} \vskip0.5cm

The next result is a version of Brezis-Lieb's Lemma
\cite{Brezis_Lieb} for Orlicz-Sobolev spaces and the proof can be
found in \cite{Gossez}.

\begin{lemma}\label{brezislieb} Let $\Omega\subset \mathbb{R}^N$ open set and $\Phi:\mathbb{R}\rightarrow [0,\infty)$
an N-function satisfies $\Delta_2-$condition. If the
complementary function $\widetilde{\Phi}$ satisfies
$\Delta_2-$condition, $(f_n)$ is bounded in $L_\Phi(\Omega)$, such
that
$$
f_n(x)\rightarrow f(x)\  \text{a.s }x\in\Omega,
$$
then
$$
f_n\rightharpoonup f \ \text{in }L_\Phi(\Omega).
$$
\end{lemma}


\section{Technical results on locally Lipschitz function and variational framework}\label{Section results for the
discontinuous}

In this section, for the reader's convenience,  we recall some
definitions and basic results  on the critical point theory of
locally Lipschitz continuous functionals  as developed by Chang
\cite{chang}, Clarke \cite{clarke, Clarke} and Grossinho \& Tersin
\cite{grossinho}. \vspace{2mm}

 Let $X$ be a real Banach space. A functional $J:X \rightarrow {\mathbb{R}}$ is locally Lipschitz continuous, $J \in
Lip_{loc}(X, {\mathbb{R}})$ for short, if given $u \in X$ there is
an open neighborhood $V := V_u \subset X$ and some constant  $K =
K_V > 0$ such that
$$
\mid J(v_2) - J(v_1) \mid \leq K \parallel v_2-v_1 \parallel,~ v_i
\in V,~ i = 1,2.
$$

\nd The directional derivative of $J$ at $u$ in the direction of $v
\in X$ is defined by

$$
J^0(u;v)=\displaystyle \displaystyle \limsup_{h \to 0,~\sigma
\downarrow 0} \frac{J(u+h+ \sigma v)-I(u+h)}{\sigma}.
$$
The generalized gradient of $J$ at $u$ is the set
$$
\partial J(u)=\big\{\mu\in X^*; \langle \mu,v\rangle\leq J^0(u;v), \ v\in
X\big \}.
$$
\nd  Since $J^0(u;0) = 0$, $\partial J(u)$ is the subdifferential of
$J^0(u;0)$. Moreover, $J^0(u;v)$ is the support function of
$\partial J(u)$ because
$$
J^0(u;v)=\max\{\langle
\xi,v\rangle; \xi\in \partial J(u)\}.
$$

The generalized gradient $\partial J(u)\subset X^*$ is convex,
non-empty and weak*-compact, and
$$
m_{J}(u) = \min \big\{\parallel\mu\parallel_{X^*};\mu \in
\partial J(u) \big \}.
$$
Moreover,
$$
\partial J(u) = \big \{J'(u) \big \},  \mbox{if}\ J \in C^1(X,{\mathbb{R}}).
$$
A critical point of $J$ is an element $u_0 \in X$ such that $0\in
\partial J(u_0)$ and a critical value of $J$ is a real number $c$
such that $J(u_0)=c$ for some critical point $u_0 \in X$.
\vspace{2mm}

About variational framework, we say that $u \in
W^{1}_{0}L_{\Phi}(\Omega)$ is a weak solution of the problem $(P_\alpha)$
if it verifies the hemivariational inequality
$$
M\left(\dis\int_\Omega\Phi(\mid\nabla u\mid )\ dx
\right)\dis\int_\Omega\phi(\mid\nabla u\mid)\nabla u\nabla v \
dx-\int_\Omega F^0(x,u;v) \ dx- \alpha\int_\Omega h(u)vdx\leq 0,
$$
for all $v \in W^{1}_{0}L_{\Phi}(\Omega)$ and moreover the set $\{x\in \Omega; \mid u\mid\geq a_0\}$ has positive measure.
Thus, weak solutions of $(P_\alpha)$ are critical points of the functional
$$
J_\alpha (u) =\widehat{M}\left(\ds\int_{\Omega}\Phi(|\nabla u|) dx\right)
-\int_{\Omega}F(x,u)dx-\alpha\int_\Omega H(u)dx,
$$
where $\widehat{M}(t)=\displaystyle\int^{t}_{0}M(s) ds$. In order to
use variational methods, we first derive some results related to the
Palais-Smale compactness condition for the problem $(P_\alpha)$.

We say that a sequence $(u_{n})\subset W^{1}_{0}L_{\Phi}(\Omega)$ is
a Palais-Smale sequence for the locally lipschitz functional $J_\alpha$
associated of problem $(P_\alpha)$ if
\begin{eqnarray}\label{****}
J_\alpha (u_{n})\rightarrow c_\alpha \ \mbox{and} \ m_{J_\alpha}(u_{n})\rightarrow 0 \
\mbox{in} \ (W^{1}_{0}L_{\Phi}(\Omega))^*.
\end{eqnarray}

If (\ref{****}) implies the existence of a subsequence $(u_{n_{j}})
\subset (u_{n})$ which converges in $W^{1,\Phi}_{0}(\Omega)$, we
say that these one functionals satisfies the nonsmooth $(PS)_{c_\alpha}$
condition.

Note that $J_\alpha \in Lip_{loc}(W^{1,\Phi}_{0}(\Omega), {\mathbb{R}})$
and from convex analysis theory, for all $ w \in \partial J_\alpha(u)$,
\begin{equation}\label{deriv}
\langle w,v\rangle=M\left(\ds\int_{\Omega}\Phi(|\nabla u|) \
dx\right) \ds\int_{\Omega}\phi(|\nabla u|)\nabla u \nabla v \ dx  -
         \langle\rho,v\rangle-\alpha\int_\Omega h(u)vdx,
\end{equation}
for some $\rho \in \partial J_F(u)$, where
$$
J_F(u)=\int_{\Omega}F(x,u) dx, \ u\in W^{1,\Phi}_{0}(\Omega).
$$
We have $J_F \in Lip_{loc}(L_{\Phi_*}(\Omega))$, $\partial J_F(u) \subset
L_{\widetilde{\Phi}_*}(\Omega)$.

The next lemma is similar to \cite[theorem~4.1]{Abrantes} and its proof will be omitted.

\begin{lemma}\label{inclusion}
Assume that $\underline{f},\overline{f}$ are N-measurable and that $(\phi_1)-(\phi_2)$ and $(f_1)$ hold. Then
\begin{description}
  \item[(i)] $\partial J_F(u)\subset \partial F(x,u), \ \text{a.e. }x\in \Omega, \ u\in L_{\Phi_*}(\Omega).$
  \item[(ii)] $\partial J_F|_{W_0^{1,\Phi}}(\Omega)(u)=\partial J_F(u), \ u\in W_0^{1,\Phi}(\Omega). $

\end{description}

\end{lemma}


\section{Results involving genus}\label{Section Genus}

We will start by considering some basic notions on the Krasnoselskii
genus that we will use in the proof of our main results.

Let $E$ be a real Banach space. Let us denote by $\mathfrak{A}$ the
class of all closed subsets  $A\subset E\setminus \{0\}$ that are
symmetric with respect to the origin, that is, $u\in A$ implies
$-u\in A$.

\begin{definition}
Let $A\in \mathfrak{A}$. The Krasnoselskii genus $\gamma(A)$ of $A$
is defined as being the least positive integer $k$ such that there
is an odd mapping $\phi \in C(A,\mathbb{R}^{k})$ such that
$\phi(x)\neq 0$ for all $x\in A$. If $k$ does not exist we set
$\gamma(A)=\infty$. Furthermore, by definition,
$\gamma(\emptyset)=0$.
\end{definition}

In the sequel we will establish only the properties of the genus
that will be used through this work. More information on this
subject may be found in the references  by  \cite{Ambrosetti},
\cite{Castro}, \cite{Clark}, \cite{Davi} and \cite{Kranolseskii}.
\begin{prop}
Let $E={\mathbb{R}}^{N}$ and $\partial\Omega$ be the boundary of an
open, symmetric and bounded subset $\Omega \subset {\mathbb{R}}^{N}$
with $0 \in \Omega$. Then $\gamma(\partial\Omega)=N$.
\end{prop}

\begin{coro}\label{esfera}
$\gamma(\mathcal{S}^{N-1})=N$ where $\mathcal{S}^{N-1}$ is a unit
sphere of ${\mathbb{R}}^{N}$.
\end{coro}

\begin{prop}\label{paracompletar}
If $K \in \mathfrak{A}$, $0 \notin K$ and $\gamma(K) \geq 2$, then
$K$ has infinitely many points.
\end{prop}


\section{Proof of Theorem \ref{teorema1}}\label{finalissimo}


The plan of the proof is to show that the set of critical points of
the functional $J_\alpha$ is compact, symmetric, does not contain the zero
and has genus  more than $2$. Thus, our main result is a consequence
of Proposition \ref{paracompletar}.

Using $(M_1)$ and Orlicz's embedding, we get
\begin{eqnarray} \label{emb 51}
J_\alpha(u)&\geq & \sigma\int_\Omega \Phi(|\nabla u|)dx- c_1\alpha
\eta_1\circ\xi_0^{-1}\left(\int_\Omega \Phi(|\nabla
u|)dx\right)\\
&&-c_2\xi_3\circ \xi_0^{-1} \left(\int_\Omega \Phi(|\nabla
u|)dx\right)\nonumber
:=g\left(\int_\Omega \Phi(|\nabla u|)dx\right),
\end{eqnarray}
where
$$
g(t)=\left \{ \begin{array}{l}
\sigma t-c_1\alpha t^{\frac{h_0}{m}}-c_2 t^{\frac{l^*}{m}}, \ \text{se } t\in (0,1]\\
\sigma t-c_1\alpha t^{\frac{h_1}{l}}-c_2 t^{\frac{m^*}{l}}, \ \text{se } t\in [1,+\infty).
\end{array}
\right.
$$

Thus, there exists $\alpha_{\#}>0$ such that, if $\alpha\in (0, \alpha_{\#})$, then $g$ attains its positive maximum. Let us assume $\alpha\in (0,\alpha_{\#})$, choosing $A_\alpha$ and $B_\alpha$ as the first and second root of $g$, we make the following of the truncation $J_\alpha$. Fix $\Psi\in C_0^\infty([0,+\infty))$,  $0\leq \Psi(t)\leq 1$ for all $t\in [0,+\infty)$, such that $\Psi(t)=1$ if $t\in [0,A_\alpha]$ and $\Psi(t)=0$ if $t\in [B_\alpha,+\infty)$. Now, we consider the truncated functional
$$
I_\alpha(u)=\widehat{M}\left(\int_\Omega \Phi(|\nabla
u|)dx\right)-\alpha\int_\Omega H(u)dx -\Psi\left(\int_\Omega
\Phi(|\nabla u|)dx\right)\int_\Omega F(x,u)dx.
$$
From $(M_1)$ and Orlicz's embedding, we have that
\begin{eqnarray*}
I_\alpha(u)&\geq& \sigma\int_\Omega \Phi(|\nabla u|)dx-\alpha\widetilde{c}_1 \eta_1\circ\xi_0^{-1}\left(\int_\Omega \Phi(|\nabla u|)dx\right)\\
&&-\widetilde{c}_2\Psi\left(\int_\Omega \Phi (|\nabla u|)dx\right)\xi_3\circ \xi_0^{-1}\left(\int_\Omega \Phi(|\nabla u|)dx\right)\\
&:=&\overline{g}\left(\int_\Omega \Phi (|\nabla u|)dx\right),
\end{eqnarray*}
where
$$
\overline{g}(t)=\left \{ \begin{array}{l}
\sigma t-\widetilde{c}_1\alpha t^{\frac{h_0}{m}}-\widetilde{c}_2 \Psi(t)t^{\frac{l^*}{m}}, \ \text{se } t\in (0,1]\\
\sigma t-\widetilde{c}_1\alpha t^{\frac{h_1}{l}}-\widetilde{c}_2 \Psi(t)t^{\frac{m^*}{l}}, \ \text{se } t\in [1,+\infty).
\end{array}
\right.
$$

Now, we will show that $I_\alpha$  satisfy the local Palais-Smale condition. For this, we need the following technical result involving the functional $J_\alpha$.

\begin{lemma}\label{PSc}
Let $(u_n)\subset W_0^{1,\Phi}(\Omega)$ be a bounded sequence such that
$$
J_\alpha(u_n)\to c_\alpha \text{ and } m^{J_\alpha}(u_n)\to 0.
$$
There exists $\alpha>0$ small enough such that, if
$$
c_\alpha<\omega_N:=\frac{1}{2}\min\left\{\left(\frac{l\sigma}{b_1 m^*}\right)^{\frac{\beta}{\beta-1}}S_N^{-\frac{\alpha}{\beta-1}}; \  \alpha\in\left\{l^*,m^*\right\} \text{ and }\beta\in
\left\{\frac{l^*}{l},\frac{m^*}{l},\frac{l^*}{m}, \frac{m^*}{m}\right\} \right\},
$$
 then up a subsequence,  $(u_n)$ is strongly in $L_{\Phi_*}(\Omega)$.
\end{lemma}
\noindent \textbf{Proof:} Since $(u_n)$ is bounded in $W_0^1L_\Phi(\Omega)$,  we can extract a subsequence of $(u_n)$, still denoted by $ (u_{n})$ and $u \in W_0^{1,\Phi}(\Omega)$, such that
\begin{itemize}
  \item $ u_{n} \rightharpoonup u $ in $ W_{0}^{1, \Phi}(\Omega) $;
  \item $ u_{n} \rightharpoonup u $ in $ L_{\Phi_*}(\Omega) $;
  \item $ u_{n} \to u $ in $ L_{\Phi}(\Omega) $;
  \item $u_n(x) \to u(x)$ a.e. in $\Omega$;
\end{itemize}

From the concentration compactness lemma of Lions in Orlicz-Sobolev space found in \cite{fukagai}, there exist two nonnegative measures
$ \mu, \nu \in \mathcal{M}(\Omega) $, a countable set $\mathcal{J}$, points $ \{ x_{j} \}_{j \in \mathcal{J}} $ in $ \Omega $ and sequences
$ \{ \mu_{j}\}_{j \in \mathcal{J}}, \{ \nu _{j}\}_{j \in \mathcal{J}} \subset [0, +\infty)$, such that
\begin{eqnarray*}
\Phi(|\nabla u_{n}|) \to \mu \geq \Phi(|\nabla u|) + \sum_{j \in \mathcal{J}}\mu_{j}\delta_{x_{j}} \mbox{ in } \mathcal{M}( \Omega)\\
\Phi_*(u_{n}) \to \nu =  \Phi_*(u) + \sum_{j \in \mathcal{J}}\nu_{j}\delta_{x_{j}} \mbox{ in } \mathcal{M}( \Omega)
\end{eqnarray*}
and
$$
\nu_j\leq\max\{S_N^{l^*}\mu_j^{\frac{l^*}{l}},S_N^{m^*}\mu_j^{\frac{m^*}{l}},S_N^{l^*}\mu_j^{\frac{l^*}{m}},
S_N^{m^*}\mu_j^{\frac{m^*}{m}}\},
$$
where $\delta_{x_j}$ is the Dirac mass at $x_j\in \Omega$ and $S_N$ verifies (\ref{trudinger-emb}).

Now, we show that there exists $\alpha_*>0$ such that
$$
\Gamma:=\{j\in \mathcal{J}; \nu_j>0\}=\emptyset, \ \text{for every } \alpha\in
(0,\alpha_*) \text{ and } c_\alpha<0.
$$
Let $ \psi \in C^{\infty}_{0}(\mathbb{R}^{N}) $ such that
$$
 \psi (x)= 1  \,\,\, \mbox{in} \,\,\,  B_{\frac{1}{2}}(0) ,  \,\,\, supp\,\psi \subset B_{1}(0) \,\, \mbox{and} \,\,\,  0 \leq \psi(x) \leq 1 \,\,\, \forall x \in \mathbb{R}^{N}.
$$
 For each $j\in \Gamma$ and $ \epsilon > 0 $, let us define
$$
\psi_{\epsilon}(x) = \psi\left(\frac{x-x_j}{\epsilon} \right),
\,\,\, \forall  x \in \mathbb{R}^{N}.
$$
Note that $ \left( \psi_{\epsilon} u_{n} \right) $ is a bounded
sequence in $W_0^{1,\Phi}(\Omega)$. Once that
$$
m_{J_\alpha}(u_n)\to 0,
$$
by (\ref{deriv})
\begin{eqnarray} \label{eqM1}
o_n(1)&=&M\left(\int_\Omega \Phi(|\nabla
u_n|)dx\right)\int_\Omega\phi(|\nabla u_n|)\nabla
u_n\nabla(u_n\psi_\epsilon)dx -\alpha\int_\Omega h(u_n)u_n\psi_\epsilon
dx\nonumber\\
&&-\langle\rho_n,\psi_\epsilon u_n\rangle,
\end{eqnarray}
where $\rho_n\in \partial J_F(u_n)$.

From (\ref{partial1}) and Lemma \ref{inclusion}
$$
\rho_n\in [\underline{f}(x,u_n),\overline{f}(x,u_n)] \ \text{a.e. } x\in \Omega,
$$
and so, by $(f_1)$
$$
|\rho_n(x)|\leq b_1 \phi_*(|u_n(x)|)|u_n(x)|,  \text{ a.e. }x\in\Omega.
$$
Using the Lemma \ref{F3}, we have that
\begin{equation}\label{eq002}
\left|\int_\Omega \rho_n u_n \psi_\epsilon dx\right|\leq b_1m^*\int_\Omega\Phi_*(u_n)\psi_\epsilon dx.
\end{equation}
Once that $(u_n)$ converges strongly for $u$ in $L_H(\Omega)$, we have
\begin{equation}\label{eq003}
\int_\Omega  H(u_n) \psi_\epsilon dx\to \int_\Omega
H(u)\psi_\epsilon dx.
\end{equation}
From (\ref{eqM1})-(\ref{eq003})
\begin{eqnarray}\label{eq004}
M\left(\int_\Omega \Phi(|\nabla u_n|)dx\right)\int_\Omega \phi(|\nabla
u_n|)\nabla u_n\nabla(u_n\psi_\epsilon)dx\leq h_1\alpha \int_\Omega
H(u)\psi_\epsilon dx\\
+b_1m_*\int_\Omega\Phi_*(u_n)\psi_\epsilon dx\nonumber+o_n(1).
\end{eqnarray}

On the other hand, by $(\phi_2)$ and $(M_1)$
\begin{align}\label{eq005}
M\left(\int_\Omega \Phi(|\nabla u_n|)dx\right)& \int_{\Omega}\phi(
|\nabla u_{n}|)\nabla u_n\nabla\left(u_n \psi_{\epsilon}\right)dx\nonumber\\
=&M\left(\int_\Omega \Phi(|\nabla u_n|)dx\right)\int_{\Omega}\phi( |\nabla u_{n}|)|\nabla u_n|^2 \psi_{\epsilon}dx\nonumber\\
&+ M\left(\int_\Omega \Phi(|\nabla u_n|)dx\right)\int_{\Omega}\phi( |\nabla u_{n}|)\left(\nabla u_n\nabla \psi_{\epsilon}\right)u_n dx\nonumber\\
\geq & l \sigma\int_\Omega \Phi(|\nabla u_n|) \psi_{\epsilon}dx\nonumber\\
&+ M\left(\int_\Omega \Phi(|\nabla u_n|)dx\right)\int_{\Omega}\phi(
|\nabla u_{n}|)\left(\nabla u_n\nabla \psi_{\epsilon}\right)u_n dx.
\end{align}
By Lemmas \ref{DESIGUALD} and \ref{lem Phiest}, the sequence
$$
\left\{\left|M\left(\int_\Omega \Phi(|\nabla u_n|)dx\right)\phi(|\nabla u_n|)\nabla
u_n\right|_{\widetilde{\Phi}}\right\}
$$ is bounded. Thus, there is a subsequence
$(u_n)$ such that
$$
M\left(\int_\Omega \Phi(|\nabla u_n|)dx\right)\phi(|\nabla u_n|)\nabla u_n\rightharpoonup \widetilde{w}_1 \mbox{
weakly in }L_{\widetilde{\Phi}}(\Omega,\mathbb{R}^N),
$$
for some $\widetilde{w}_1\in
L_{\widetilde{\Phi}}(\Omega,\mathbb{R}^N)$. Since  $u_n \to u$ in
$L_\Phi( \Omega)$,
$$
M\left(\int_\Omega \Phi(|\nabla u_n|)dx\right)\int_\Omega \phi(|\nabla u_n|)(\nabla u_n\nabla \psi_\epsilon)u_n dx\to
\int_\Omega (\widetilde{w}_1\nabla\psi_\epsilon)u dx.
$$
Thus, combining (\ref{eq004}), (\ref{eq005}) and letting $n\to
\infty$, we have
\begin{equation}\label{eq006}
l\sigma\int_\Omega \psi_\epsilon d\mu+\int_\Omega
(\widetilde{w}_1\nabla \psi_\epsilon)u dx \leq \alpha h_1 \int_\Omega
H(u)\psi_\epsilon dx  +b_1m^*\int_\Omega \psi_\epsilon d\nu.
\end{equation}

Now we show that
$$
\int_\Omega (\widetilde{w}_1\nabla \psi_\epsilon)u dx\to 0,
\text{as }\epsilon \to 0.
$$
Once that $|\{\alpha
h(u_n)u_n+\rho_n\}|_{\widetilde{\Phi}_*}$ is bounded, let
$\widetilde{w}_2\in L_{\widetilde{\Phi}_*}(\Omega)$ such that
$$
\alpha h(u_n)u_n+\rho_n\rightharpoonup \widetilde{w}_2 \text{ weakly in
}L_{\widetilde{\Phi}_*}(\Omega).
$$
Since
$$
M\left(\int_\Omega \Phi(|\nabla u_n|)dx\right)\int_\Omega
\phi(|\nabla u_n|)\nabla u_n\nabla v dx-\int_\Omega (\alpha
h(u_n)u_n+\rho_n)v dx\to 0,
$$
as $n\to+\infty $ for any $v\in W_0^{1,\Phi}(\Omega)$,
$$
\int_\Omega \widetilde{w}_1\nabla v- \widetilde{w}_2v dx=0,
$$
for any $v\in W_0^{1,\Phi}(\Omega)$. Substituting $v=u\psi_\epsilon$ we have
$$
\int_\Omega \widetilde{w}_1\nabla (u\psi_\epsilon)-
\widetilde{w}_2u\psi_\epsilon dx=0.
$$
Namely,
$$
\int_\Omega (\widetilde{w}_1\nabla \psi_\epsilon)u dx=-\int_\Omega
(\widetilde{w}_1\nabla u-\widetilde{w}_2u)\psi_\epsilon dx.
$$
Noting $\widetilde{w}_1\nabla u-\widetilde{w}_2 u\in L^1(\Omega)$, we see that right-hand side tends to 0 as $\epsilon\to 0$. Hence we have
$$
\int_\Omega \left(\widetilde{w}_1\nabla \psi_\epsilon\right) u dx \to 0,
$$
as $\epsilon \to 0$.
Letting $\epsilon \to 0$ in (\ref{eq006}), we obtain
$$
l\sigma\mu_j\leq b_1m^*\nu_j.
$$
Hence
$$
S_N^{-\alpha}\nu_j\leq \mu_j^\beta \leq
\left(\frac{b_1m^*}{l\sigma}\right)^\beta \nu_j^{\beta},
$$
for some $\alpha\in\{l^*,m^*\}$ and $\beta\in
\{\frac{l^*}{l},\frac{m^*}{l},\frac{l^*}{m}, \frac{m^*}{m}\}$, and so
\begin{equation}\label{eq008}
\nu_j\geq
\left(\frac{l\sigma}{b_1 m^*}\right)^{\frac{\beta}{\beta-1}}S_N^{-\frac{\alpha}{\beta-1}}>0
\text{ or }\nu_j=0.
\end{equation}

Now, taking $\theta\in(2m,l^*)$, we get
\begin{eqnarray*}
c_\alpha&=&\displaystyle \lim_{n\to
+\infty}\left(J_\alpha(u_n)-\frac{1}{\theta}J'_\alpha(u_n)u_n\right)\\
&&\geq \displaystyle \lim_{n\to
+\infty}\left(\widehat{M}\left(\int_\Omega \Phi(|\nabla
u_n|)dx\right)-\frac{m}{\theta}M\left(\int_\Omega \Phi(|\nabla
u_n|)dx\right)\right)\int_\Omega \Phi(|\nabla u_n|)dx\\
&&-\alpha\left(1-\frac{h_0}{\theta}\right)\int_\Omega H(u_n)dx-\int_\Omega\left(
F(x,u_n)-\frac{1}{\theta}\rho_nu_n \right)dx,
\end{eqnarray*}
and so by  $(M_2)$ and $(f_1)$
\begin{eqnarray*}
c_\alpha&\geq&\displaystyle \lim_{n\to +\infty}\left(-\alpha
\left(1-\frac{h_0}{\theta}\right)\int_\Omega
H(u_n)dx+\left(\frac{b_0l^*}{\theta}-b_1\right)\int_\Omega\Phi_*(u_n)\right),\\
&&=-\alpha \left(1-\frac{h_0}{\theta}\right)\int_\Omega H(u)dx+\left(\frac{b_0l^*}{\theta}-b_1\right)\left(\int_\Omega \Phi_*(u)dx+\nu_j\right),
\end{eqnarray*}
where $\frac{b_0l^*}{\theta}-b_1>0 $ (see $(f_1)$). Consequently
\begin{equation}\label{eq009}
\omega_N>c_\alpha\geq -\alpha \left(1-\frac{h_0}{\theta}\right)\eta_1(|u|_{\Phi_*})+\left(\frac{b_0l^*}{\theta}-b_1\right) \left(\xi_3(|u|_{\Phi_*})+
\nu_j\right).
\end{equation}
Define
$$
\zeta(t)=-\alpha \left(1-\frac{h_0}{\theta}\right)\eta_1(t)+\left(\frac{b_0l^*}{\theta}-b_1\right) \xi_3(t), \ t\geq 0,
$$
note that
\begin{equation*}
\zeta(t)=\left\{
\begin{array}{l}
-\alpha \left(1-\frac{h_0}{\theta}\right)t^{h_0}+\left(\frac{b_0l^*}{\theta}-b_1\right)t^{m^*}, \ \text{ se
 }t\in(0,1]\\
-\alpha \left(1-\frac{h_0}{\theta}\right)t^{h_1}+\left(\frac{b_0l^*}{\theta}-b_1\right)t^{l^*}, \ \text{ se }t\in[1,+\infty).
\end{array}
\right.
\end{equation*}
This function attains its absolute minimum, for $t>0$, at the point
$$
t_0:=\left(\frac{\alpha \left(\theta-h_0\right)
h_0}{\left(b_0l^*-b_1\theta\right)m^*}\right)^{\frac{h_0-1}{m^*-1}}.
$$
Thus, we conclude for $\alpha>0$
small enough that
\begin{equation}\label{eq010}
\omega_N+\alpha^{h_0\left(\frac{h_0-1}{m^*-1}\right)+1}\left(1-\frac{h_0}{\theta}\right)\lambda_0^{h_0}
-\alpha^{m^*\left(\frac{b_0-1}{m^*-1}\right)}\left(\frac{b_0l^*}{\theta}-b_1\right)\lambda_0^{m^*}
>\nu_j,
\end{equation}
where $\lambda_0:=\left(\frac{\left(\theta-h_0\right)
h_0}{\left(b_0l^*-b_1\theta\right)m^*}\right)^{\frac{h_0-1}{m^*-1}}>0.$

Supposing that there exists $j\in \Gamma$ such that $\nu_j>0$. From (\ref{eq008}) and (\ref{eq010})
$$
\left(1-\frac{h_0}{\theta}\right)\lambda_0^{h_0}\alpha^{h_0\left(\frac{h_0-1}{m^*-1}\right)+1}
-\left(\frac{b_0l^*}{\theta}-b_1\right)\lambda_0^{m^*}\alpha^{m^*\left(\frac{h_0-1}{m^*-1}\right)}>\frac{1}{2}\omega_N,
$$
but this is a contradiction, for $\alpha>0$ small enough.
Showing that
$$
u_n \to u \text{ in }L_{\Phi_*}(\Omega).
$$
\hfill\rule{2mm}{2mm}
\begin{lemma}\label{LPS}
Let $(u_n)\subset W_0^{1,\Phi}(\Omega)$ the $(PS)_{c_\alpha}$ sequence obtained in the previous. Then, for some subsequence, still denoted by it self,
$$
u_n\to u \text{ in }W_0^{1,\Phi}(\Omega).
$$
\end{lemma}
 \textbf{Proof:} Now, as $ m_{J_\alpha}(u_n) = o_{n}(1) $, the last limit gives
\begin{align*}
M\left(\int_\Omega \Phi(|\nabla
u_n|)dx\right)\int_{\Omega}\phi(|\nabla u_n|)\mid u_n\mid^2dx  =& \alpha
\int_{\Omega}h(u_n)u_n^2dx  + \int_{\Omega}\rho_n u_{n}dx+ o_n(1).
\end{align*}
In what follows, let us denote by $\{P_n\}$ the following sequence,
$$
P_{n}(x) = \langle \phi(|\nabla u_{n}(x)|) \nabla u_{n}(x) -
\phi(|\nabla u(x)|) \nabla u(x), \nabla u_{n}(x) - \nabla u (x)
\rangle.
$$
Since $\Phi$ is convex in $\mathbb{R}$ and $\Phi(|.|)$ is $C^1$
class in $\mathbb{R}^N$,  has $P_n(x) \geq 0. $ From definition of
$\{P_n\}$,
\[
\int_{\Omega} P_{n}dx = \int_{\Omega} \phi(|\nabla u_{n}|)\mid\nabla
u_n\mid^2dx - \int_{\Omega}\phi( |\nabla u_{n}|) \nabla u_{n} \nabla udx
- \int_{\Omega} \phi(|\nabla u|) \nabla u \nabla(u_{n} - u)dx.
\]
Recalling that $u_n \rightharpoonup u$ in $W_0^{1,\Phi}(\Omega)$, we
have
\begin{align}
\int_{\Omega} \phi(|\nabla u|) \nabla u \nabla(u_{n} - u) dx \to 0
\quad \mbox{ as } n \to \infty,
\end{align}
which implies that
\[
\int_{\Omega} P_{n}dx = \int_{\Omega} \phi(|\nabla u_{n}|) |\nabla
u_{n}|^2dx- \int_{\Omega} \phi(|\nabla u_{n}|)\nabla u_{n} \nabla udx +
o_n(1).
\]
On the other hand
\begin{align*}
0\leq M\left(\int_\Omega \Phi(|\nabla u_n|)dx\right)\int_{\Omega} P_{n}dx =& \alpha \int_{\Omega} h(u_{n})u_n^2dx  - \alpha \int_{\Omega}h(u_n)u_{n} udx\\
    &  +\int_{\Omega} \rho_n u_{n}dx - \int_{\Omega}\rho_n udx + o_n(1),
\end{align*}
where $\rho_n\in\partial J_F(u_n)$. Once that $|\rho_n|_{\widetilde{\Phi}_*}$ is bounded, we have that
\begin{align*}
\int_{\Omega} P_{n}dx&\leq \frac{\alpha}{\sigma}\int_\Omega
h(u_n)u_n(u_n-u)dx+\frac{1}{\sigma}\int_\Omega \rho_n(u_n-u)dx\\
&\leq
\frac{\alpha}{\sigma}|u_n|_{\widetilde{H}}|u_n-u|_{H}+\frac{1}{\sigma}|\rho_n|_{\widetilde{\Phi}_*}|u_n-u|_{\Phi_*}\to
0.
\end{align*}
 Applying a result due to Dal Maso and Murat
\cite{Maso}
$$
u_{n} \to u \mbox{ in } W^{1,\Phi}_{0}(\Omega) .
$$
\hfill\rule{2mm}{2mm}

\begin{lemma}\label{igualdade}
If $I_\alpha(u)<0$, then $\int_\Omega\Phi(|\nabla u|)dx<A_\alpha$ and $J_\alpha(v)=I_\alpha(v)$, for all $v$ in a sufficiently small neighborhood of $u$. Moreover, $I_\alpha$ verifies a local Palais-Smale condition for $c_\alpha<0$.
\end{lemma}
 \textbf{Proof:} Since $\overline{g}(\int_\Omega\Phi(|\nabla u|)dx)\leq I_\alpha(u)<0$, we have that  $\int_\Omega\Phi(|\nabla u|)dx<A_\alpha$. Once that $\{v\in W_0^{1,\Phi}(\Omega);  \int_\Omega\Phi(|\nabla v|)dx<A_\alpha\}$ is open, there exists $r_\alpha>0$ such that  $\int_\Omega\Phi(|\nabla v|)dx<A_\alpha$ for all $v\in B_{r_\alpha}(u)$ arguing as in Theorem \ref{teorema1} we conclude that $J_\alpha(v)=I_\alpha(v)$, for all $v\in B_{r_\alpha}(u)$. Moreover, if $(u_n)$ is a sequence such that $I_\alpha(u_n)\to c_\alpha<0$ and $m_{I_\alpha}(u_n)\to 0$, then, for $n$ sufficiently large, $J_\alpha(u_n)=I_\alpha(u_n)\to c_\alpha<0$ and $m_{J_\alpha}(u_n)=m_{I_{\alpha}}(u_n)\to 0$. Since $I_\alpha$ is coercive, we get that $(u_n)$ is bounded in $W_0^{1,\Phi}(\Omega)$. From Lemma \ref{LPS}, for $\alpha$ sufficiently small and $c_\alpha<0$, hence, up to a subsequence, $(u_n)$ is strongly convergent in $W_0^{1,\Phi}(\Omega)$.
\hfill\rule{2mm}{2mm}

Let $K_{c}$ be the set of critical points of $J$. More precisely
$$
K_{c}(J)=\{u \in W^{1}_{0}L_{\Phi}(\Omega): 0 \in \partial J(u) \ \
\mbox{and} \ \ J(u)=c\}.
$$

Since $I_\alpha$ is even, we have that $K_{c_\alpha}:=K_{c_\alpha}(I_\alpha)$ is symmetric. The next result is important in our arguments and allows  we conclude that
$K_{c_\alpha}$ is compact. The proof can be found in \cite{chang}.
\begin{lemma}\label{KcCompacto}
If $I_\alpha$ satisfies the nonsmooth $(PS)_{c_\alpha}$ condition, then $K_{c_\alpha}$ is compact.
\end{lemma}

To prove that $K_{c_\alpha}$ does not contain zero, we construct a special
class of the levels $c_\alpha$.

For each $k \in \mathbb{N}$, we define the set
$$
\Gamma_{k}=\{C \subset W^{1}_{0}L_{\Phi}(\Omega): C \ \ \mbox{is
closed}, C=-C \ \ \mbox{and} \ \ \gamma(C) \geq k\},
$$
and the values
$$
c^{\alpha}_{k}=\displaystyle\inf_{C\in \Gamma_{k}}\displaystyle\sup_{u \in
C}J_\alpha(u).
$$

Note that
$$
-\infty\leq c^{\alpha}_{1}\leq c^{\alpha}_{2}\leq c^{\alpha}_{3}\leq ...\leq
c^{\alpha}_{k}\leq ...
$$
and, once that $I_\alpha$ is coercive and continuous, $I_\alpha$ is bounded below
and, hence, $c^{\alpha}_{1} > -\infty$. In this case, arguing as in
\cite[Proposition 3.1]{BWW}, we can prove that each $c^{\alpha}_{k}$ is a
critical value for the functional $I_\alpha$.

\begin{lemma}\label{minimax}
Given $k \in \mathbb{N}$, there exists $\epsilon = \epsilon(k)>0$
such that
$$
\gamma(I_\alpha^{-\epsilon}) \geq k,
$$
where $I_\alpha^{-\epsilon}=\{u \in W^{1}_{0}L_{\Phi}(\Omega): I_\alpha(u) \leq
-\epsilon\}$.
\end{lemma}
\noindent\textbf{Proof:} Fix $k \in \mathbb{N}$, let $X_{k}$ be a
k-dimensional subspace of $W^{1}_{0}L_{\Phi}(\Omega)$. Thus, there
exists $C_k>0$ such that
$$
-C_k\parallel u\parallel_{\Phi}\geq - \displaystyle|u|_{H},
$$
for all $u \in X_{k}$.

Thus, for all $R\in(0,1)$ and for all $u\in X_k$ with $\int_\Omega\Phi(|\nabla u|)dx<R$, from continuity of the function $M$, we conclude  that, there exists $C$, such that
\begin{eqnarray*}
I_\alpha(u)&\leq& C \xi_1(\parallel u\parallel_\Phi)-\alpha\eta_0(C_k\parallel u\parallel_\Phi)\\\nonumber
&=&C\parallel u\parallel_\Phi^l-\alpha C_k^{h_1}\parallel u\parallel_\Phi^{h_1}.
\end{eqnarray*}

Fixing $r\in (0,1)$ such that
$$
C r^{l-h_1}<\alpha C_k^{h_1},
$$
and $ {\mathcal{S}_r}=\{u\in X_k; \parallel u \parallel_\Phi=r \}$, we get
\begin{eqnarray*}
I_\alpha(u)&\leq& C r^l-\alpha C_k r^{h_1}<0=I_\alpha(0), \ \forall u\in {\mathcal{S}_r},
\end{eqnarray*}
which implies there exists $\epsilon=\epsilon(r)>0$ such that
$$
I_\alpha(u)<-\epsilon < 0,
$$
for all $u\in {\mathcal{S}_r}$.
Since $X_k$ and $\mathbb{R}^k$ are isomorphic and $\mathcal{S}_r$
and $S^{k-1}$ are homeomorphic, we conclude from Corollary
\ref{esfera} that $\gamma(\mathcal{S}_r)=\gamma(S^{k-1})=k$.
Moreover, once that ${\mathcal{S}_r} \subset I_\alpha^{-\epsilon}$ and
$I_\alpha^{-\epsilon}$ is symmetric and closed,  we have
$$
k= \gamma ({\mathcal{S}_r})\leq \gamma( I_\alpha^{-\epsilon}).
$$
\hfill\rule{2mm}{2mm}

\begin{lemma}\label{minimax1}
Given $k \in \mathbb{N}$, the number $c^{\alpha}_{k}$ is negative.
\end{lemma}
\noindent\textbf{Proof:} From Lemma \ref{minimax}, for each $k\in
\mathbb{N}$ there exists $\epsilon >0$ such that
$\gamma(I_\alpha^{-\epsilon}) \geq k$. Moreover, $ 0 \notin I_\alpha^{-\epsilon}$
and $I_\alpha^{-\epsilon}\in \Gamma_{k}$. On the other hand
$$
\displaystyle\sup_{u\in I_\alpha^{-\epsilon}}I_\alpha(u)\leq -\epsilon.
$$

Hence,
$$
-\infty < c^{\alpha}_{k}=\displaystyle\inf_{C\in
\Gamma_{k}}\displaystyle\sup_{u \in C}I_\alpha(u) \leq
\displaystyle\sup_{u\in I_\alpha^{-\epsilon}}I_\alpha(u) \leq -\epsilon <0.
$$
\hfill\rule{2mm}{2mm}

The next result is a direct adaptation of \cite[Lemma 4.4]{GP}. See
also \cite[Lemma 7.5]{GJ}.

\begin{lemma}\label{minimax2}
If $c^{\alpha}_{k}=c^{\alpha}_{k+1}=...=c^{\alpha}_{k+r}$ for some $r \in \mathbb{N}$, then
$$
\gamma(K_{c^{\alpha}_{k}})\geq r+1.
$$
\end{lemma}

\subsection{Proof of Theorem \ref{teorema1}}

If $-\infty< c^{\alpha}_{1} < c^{\alpha}_{2} < ...< c^{\alpha}_{k}< ...<0$ and since each
$c^{\alpha}_{k}$ critical value of $I^{\alpha}=J^{\alpha}$ (see Lemma \ref{igualdade}), then we obtain infinitely many
critical points of $J_\alpha$.

On the other hand, if there are two constants $c^{\alpha}_{k}=c^{\alpha}_{k+r}$, then
$c^{\alpha}_{k}=c^{\alpha}_{k+1}=...=c^{\alpha}_{k+r}$ and from Lemma \ref{minimax2}, we have
$$
\gamma(K_{c^{\alpha}_{k}})\geq r+1 \geq 2.
$$
From Proposition \ref{paracompletar}, $K_{c^{\alpha}_{k}}$ has infinitely
many points.

Let $(u_{k})$ critical points of $J_\alpha$. Now we show that, for
\begin{eqnarray}\label{setmeasure}
a_0 < \eta_1^{-1}\left(\frac{\sigma l}{H(1)\mid \Omega\mid
\alpha h_1}\xi_0(C)\right),
\end{eqnarray}
we have that
$$
\bigl\{x \in \Omega: \mid u_{k}(x)\mid \geq a_0\bigl\}
$$
has positive measure. Thus every critical points of $J_\alpha$, are solutions of $(P)$. Suppose, by contradiction, that this set has
null measure. Thus
\begin{eqnarray*}
0&=&M\left(\int_\Omega\Phi(\mid \nabla u_k\mid )dx\right)\int_\Omega
\phi(\mid\nabla u_k\mid)\mid\nabla u_k\mid^2dx-\alpha\int_\Omega h(
u_k)dx \\
  &\geq &  \sigma l \int_\Omega\Phi(\mid\nabla
u_k\mid)dx-\alpha h_1\int_\Omega H(u_k)dx\\
 &\geq & \sigma l\xi_0(\parallel
u_k\parallel_\Phi)-\alpha h_1H(a_0)\mid \Omega\mid,\
 \end{eqnarray*}
where we conclude
\begin{eqnarray}\label{setmeasure1}
\sigma l\xi_0(\parallel u_k\parallel_\Phi)^{(\alpha+1)}\leq \alpha h_1\eta_1(a_0)\mid \Omega\mid H(1).
\end{eqnarray}

Since $c^{\alpha}_{k}\leq -\epsilon<0$, there exists $C>0$ such that
$\|u_{k}\|\geq C>0$. Hence
\begin{eqnarray*}
a_0 \geq \eta_1^{-1}\left(\frac{\sigma l}{H(1)\mid \Omega\mid
\alpha h_1}\xi_0(C)\right),
\end{eqnarray*}
which contradicts (\ref{setmeasure}). Then, $$ \bigl\{x \in \Omega:
\mid u_{k}(x)\mid \geq a_0\bigl\}
$$
has positive measure. \hfill\rule{2mm}{2mm}


\begin{thebibliography}{99}

\bibitem{adams} A. Adams and J. F. Fournier, Sobolev spaces, {\it 2nd ed.}, Academic Press, (2003).


\bibitem{alvescorreama} C.O. Alves, F.J.S.A. Corr\^{e}a and T.F Ma, {\it
Positive solutions for a quasilinear elliptic equation of Kirchhoff
type}, Comput. Math. Appl., 49(2005)85-93.






\bibitem{alvescorrea} C.O. Alves and  F.J.S.A. Corr\^{e}a , {\it On
existence of solutions for a class of problem involving a nonlinear
operator}, Comm. Appl. Nonlinear Anal., 8(2001)43-56.

\bibitem{Abrantes} {C.O. Alves, J.V. Gon\c calves and J.A. Santos,} {\it Strongly Nonlinear Multivalued Elliptic
Equations on a Bounded Domain}, J Glob Optim, DOI
10.1007/s10898-013-0052-3, 2013.

\bibitem{Ambrosetti} A. Ambrosetti and P. H Rabinowitz,
{\it Dual variational methods in cri\-ti\-cal point theory and
apllications}, J. Functional Analysis, vol 14(1973)349-381.

\bibitem{Arosio} A. Arosio, {\it On the nonlinear Timoshenko-Kirchoff beam equation},
Chin. Annal Math., 20 (1999), 495-506.

\bibitem{Arosio1} A. Arosio, {\it A geometrical nonlinear correction to the Timoshenko beam equation},
 Nonlinear Anal. 47(2001), 729-740.

 \bibitem{GP} J. G. Azorero and I. P. Alonso,
{\it Multiplicity of solutions for elliptic problems with critical
exponent or with a nonsymmetric term}, Trans. Amer. Math. Soc. , vol
323 n. 2(1991)877-895.

\bibitem{BWW} T. Bartsch, T. Weth and M. Willem, {\it Partial symmetry
of least energy nodal solution to some variational problems}, J.
d'Analyse Mathematique, 96(2005), 1-18.

\bibitem{Bonanno} G. Bonanno, G. M. Bisci and V. Radulescu,
{\it Quasilinear elliptic non-homogeneous Dirichlet problems through
Orlicz�Sobolev spaces}, Nonlinear Analysis 75 (2012) 4441-4456.


\bibitem{Brezis_Lieb} H. Brezis. and E. Lieb, {\it A relation between pointwise convergence of functions and
convergence of functinals}, Proc. Amer. Math. Soc. 88 (1983),
486-490.

\bibitem{Castro} A. Castro,
{\it Metodos variacionales y analisi functional no linear}, X
Coloquio colombiano de Matematicas, 1980.

\bibitem{Carvalho} M.L. Carvalho and J.V. Goncalves, {\it Multivalued
Equations on a Bounded Domain via Minimization on Orlicz-Sobolev},
Journal of Convex Analysis, preprint 2014.

\bibitem{chang}  K.C. Chang, {\it Variational methods for nondifferentiable
functionals and their applications to partial differential
equations,} J. Math. Anal., 80 (1981)102-129.

\bibitem{chang1} K. C. Chang, {\it On the multiple solutions of the elliptic differential equations with discontinuous
nonlinear terms}  Sci. Sinica 21 (1978) 139-158.

\bibitem{chang2} K. C. Chang, {\it The obstacle problem and partial differential equations with discontinuous
nonlinearities}  Comm. Pure Appl. Math (1978) 139-158.

\bibitem{Chung} N. T. Chung and H. Q. Toan, {\it On a nonlinear and
non-homogeneous problem without (A-R) type condition in
Orlicz�Sobolev spaces},  Appl. Math.  and Comput. 219 (2013)
7820-7829.

\bibitem{Chung1} N. T. Chung, {\it Multiple solutions for a nonlocal problem
in Orlicz-Sobolev spaces},  Ricerche mat (2013), DOI
10.1007/s11587-013-0171-7.

\bibitem{Clark} D.C. Clark, {\it A variant of the Lusternik-Schnirelman theory}, Indiana Univ. Math. J., 22(1972)65-74.

\bibitem{clarke} F.H. Clarke,{\it Optimization and Nonsmooth Analysis}, John Wiley \& Sons, N.Y, 1983.

\bibitem{Clarke} F.H. Clarke,{\it Generalized gradients and applications}, Trans. Amer. Math. Soc. 265 (1975), 247-262.


\bibitem{correa} F.J. S. A. Corr\^{e}a and G. M. Figueiredo, {\it On a p-Kirchhoff equation via Krasnoselskii's genus},
Applied Math. Letters, 22(2009)819-822.

\bibitem{Nascimento} F.J. S. A. Corr\^{e}a and R. G. Nascimento, {\it Existence of solutions to nonlocal elliptic equations with discontinuous terms
}, EJDE 26(2012) 1-14.


\bibitem{Davi} D. G. Costa {\it An invitation to variational methods in Differential Equations}, Birkha\"{u}ser Boston, 2007.

\bibitem{Maso}G. Dal Maso and F. Murat, {\it \,  Almost everywhere convergence of gradients of solutions
to nonlinear elliptic systems,} Nonlinear Anal. 31 (1998), 405-412.

\bibitem{Mihailescu} M. Mihailescu and D. Repovs,
{\it  Multiple solutions for a nonlinear and non-homogeneous problem
in Orlicz-Sobolev spaces,} Appl. Math. and Comput. 217 (2011)
6624-6632.



\bibitem{Donaldson2}{T.K. Donaldson and N.S. Trudinger,}{\it \, Orlicz-Sobolev spaces and imbedding theorems,} J. Funct. Anal. 8 (1971) 52-75.

\bibitem{GJ} G. M. Figueiredo and J.R. dos Santos Junior,
{\it Multiplicity of solutions for a Kirchhoff equation with
subcritical or critical growth},  DIE-Diff. Int. Equations,
25(2012), 853-868.

\bibitem{fukagai}N. Fukagai and K. Narukawa,
{\it Positive solutons of quasilinear elliptic equations with
critical Orlicz-Sobolev nonlinearity on $\mathbb{R}^{N}$},
Funkciallaj Ekvacioj, 49(2006)235-267.



\bibitem{grossinho} M. R. Grossinho and S. A. Tersian {\it
An Introduction to Minimax theorems and their Applications to
Differential Equations}, 2001.

\bibitem{Gossez} J.P. Gossez {\it Orlicz-Sobolev spaces and nonlinear elliptic boundary
value problems.} In: Fuk, Svatopluk and Kufner, Alois (eds.):
Nonlinear Analysis, Function Spaces and Applications, Proceedings of
a Spring School held in Horn Bradlo, 1978. [Vol 1]. BSB B. G.
Teubner Verlagsgesellschaft, Leipzig, 1979. pp. 59-94.



\bibitem{Xiao} {X. He and W. Zou}, {\it Multiplicity of Solutions for a Class of Kirchhoff Type Problems},
  Acta Math. Applicatae Sinica, 26(2010), 387-394.

\bibitem{Kranolseskii} M. A. Kranolseskii, {\it Topological methods in the theory of nonlinear integral equations}. MacMillan, New York, 1964.

\bibitem{kirchhoff}
G. Kirchhoff {\it Mechanik, Teubner,Leipzig, 1883}.

\bibitem{lions} J.L. Lions
{\it On some questions in boundary value problems of mathematical
physics} {\it International Symposium on Continuum, Mechanics and
Partial Differential Equations, Rio de Janeiro(1977), Mathematics
Studies, Vol. 30, North-Holland, Amsterdam (1978)284-346}.

\bibitem{ma} T. F. Ma, {\it Remarks on an elliptic equation of
Kirchhoff type}. Nonlinear Anal., 63(2005,)1967-1977.





 \bibitem{Ji} J. Sun and C. Tang, {\it Existence and multiplicity of solutions for Kirchhoff type equations,} Nonlinear
 Anal., 74 (2011), 1212-1222.

  \bibitem{Santos} J. A. Santos, {\it Multiplicity of solutions for quasilinear
equations involving critical Orlicz-Sobolev nonlinear terms,} EJDE,
Vol. 2013 (2013), No. 249, pp. 1-13.


\end{thebibliography}
\end{document}